\documentclass[10pt]{article}
\usepackage[all]{xy}
\usepackage{amsmath}
\usepackage{amsfonts}
\usepackage{amssymb}
\usepackage{amscd}
\usepackage{amsthm}
\usepackage{latexsym}
\usepackage{amsbsy}

\begin{document}

\title{\bf A Noncommutative Residue on Tori and a Semiclassical  Limit}

\author {{\bf Farzad Fathizadeh}}

\date{}
\maketitle

\begin{abstract}
We define a noncommutative residue for classical 
Euclidean pseudodifferential operators on  a torus of arbitrary dimension. We prove that, 
up to multiplication by a constant, it is the unique trace on the algebra of classical 
pseudodifferential operators modulo infinitely smoothing 
operators. In the case of the two torus, we show that the noncommutative 
residue  is the semiclassical limit of a noncommutative residue defined on classical 
pseudodifferential operators on noncommutative two tori. 
\end{abstract}



\section{Introduction}

Original works on noncommutative residues as traces on algebras of 
pseudodifferential operators can be traced back to 
\cite{Adl, Man} of Adler and Manin on one dimensional symbols. In a remarkable work, Wodzicki defined 
the noncommutative residue in higher dimensions  and proved that it is the unique trace on the algebra of 
pseudodifferential operators on compact manifolds \cite{Wod}. The Wodzicki noncommutative residue bas been 
generalized vastly in the context of the local index formula in noncommutative geometry 
by Connes and Moscovici \cite{conmos}. 
In fact, using residue trace functionals and assuming the \emph{simple discrete dimension spectrum hypothesis}, 
they show that the generalization of Wodzicki's residue is a trace on the algebra of pseudodifferential 
operators associated to a \emph{spectral triple}.

In this paper  we define a noncommutative residue 
for classical {\it Euclidean} pseudodifferential operators on tori \cite{RuzTur} and 
prove that up to a constant multiple, it is the unique {\it continuous} trace on the algebra of classical 
pseudodifferential operators modulo infinitely smoothing 
operators. We also show that for the two torus, our noncommutative residue is  the semiclassical limit of a 
noncommutative residue defined on classical 
pseudodifferential operators on noncommutative two tori \cite{FatWon}.

 The author would like to thank IHES for kind support and excellent environment during 
 his visit in Summer 2011 where part of this work was carried out.

\section{Euclidean Symbols on Tori}
In this section we briefly recall a class of pseudodifferential operators on tori. We denote by $S^m(\mathbb{T}^n\times \mathbb{R}^n)$ 
the set of \emph{Euclidean symbols} 
of order $m \in \mathbb{Z}$  on the torus $\mathbb{T}^n = 
\mathbb{R}^n / 2 \pi  \mathbb{Z}^n$ \cite{RuzTur}.  This is the set of  all $C^\infty$ 
functions $\sigma : 
\mathbb{T}^n \times \mathbb{R}^n \to \mathbb{C}$ with 
the property that for all multi-indices $\alpha$ and  
$\beta$, there exists a positive constant $C_{\alpha, \beta}$ such 
that
$$
     | (\partial_{x}^{\alpha} \partial_{\xi}^{\beta} \sigma) (x, \xi)  | 
\leq C_{\alpha, \beta} \langle  \xi \rangle ^ {m - |\beta|} $$ 
for all $(x, \xi) \in \mathbb{T}^n \times \mathbb{R}^n$, where $\langle  
\xi \rangle = (1 + |\xi|^2)^{1/2}$. 
Throughout this paper, we assume 
that $n \geq 2$. 
 
Let $\sigma\in S^m(\mathbb{T}^n\times \mathbb{R}^n)$. Then the 
corresponding \emph{
pseudodifferential operator} on $\mathbb{T}^n$ is defined by 
\begin{equation} \label{pseudodef}
(T_{\sigma} f) (x) = \int_{\mathbb{R}^n} \int_{\mathbb{T}^n} e^{ i (x-y) 
\cdot \xi} \sigma(x, \xi) f(y) \, \text{d}y \, \text{d}\xi, \quad x\in 
\mathbb{T}^n,
\end{equation} 
for all smooth functions $f$ on $\mathbb{T}^n.$
Pseudodifferential operators with symbols in $\cup_{m\in \mathbb{R}}S^m(\mathbb{T}^n\times 
\mathbb{R}^m)$ form an algebra \cite{RuzTur}. This means that if 
$\sigma \in S^{m_1}(\mathbb{T}^n  \times \mathbb{R}^n)$ and $\tau 
\in S^{m_2}(\mathbb{T}^n  \times \mathbb{R}^n)$, then 
$T_{\sigma} T_{\tau}$ is a pseudodifferential operator of order 
$m_1 + m_2$ of which the  symbol $\lambda \in 
S^{m_1+m_2}(\mathbb{T}^n  \times \mathbb{R}^n)$ 
has an  asymptotic expansion given by
\begin{equation} \label{compositionsymbol}
\lambda \sim \sum_{\gamma} \frac{1}{\gamma !} 
(\partial_{\xi}^\gamma\sigma) (D_x^{\gamma}\tau). 
\end{equation}

A symbol $\sigma \in S^m(\mathbb{T}^n \times \mathbb{R}^n)$ is said to be 
\emph{classical} if it admits an asymptotic expansion 
of the form
\begin{equation} \label{classicaldef}
\sigma(x,\xi) \sim \sum_{j=0}^{\infty} \sigma_{m-j}(x,\xi)    \,\,\, 
\textnormal{as} 
\,\,\, |\xi| \to \infty ,
\end{equation}
where each $\sigma_{m-j}: \mathbb{T}^n \times (\mathbb{R}^n \setminus 
\{0 \}) \to \mathbb{C} $ is $C^\infty$ and \emph{positively  
homogeneous} in $\xi$ of order $m-j$, \emph{i.e.,}
$$
\sigma_{m-j} (x, t \xi) = t^{m-j} \sigma_{m-j} (x, \xi) $$
for all $(x, \xi) \in \mathbb{T}^n \times (\mathbb{R}^n \setminus \{0 \})$ 
and $t > 0$. The set of 
classical symbols of order $m$ and the corresponding set of 
pseudodifferential operators are denoted by 
$S_{{\rm cl}}^m(\mathbb{T}^n \times \mathbb{R}^n)$ and 
$\Psi_{{\rm cl}}^m(\mathbb{T}^n 
\times \mathbb{R}^n)$ respectively. Using a 
similar argument to the one given in \cite{FatWon}, we can see that the 
homogeneous terms in the above asymptotic expansion 
are uniquely determined by $\sigma$. The pseudodifferential operators 
associated with these symbols are also called 
\emph{classical}.

\section{ A Noncommutative Residue} \label{residueeuclidean}
\setcounter{equation}{0}
In this section we define a noncommutative residue for classical 
Euclidean pseudodifferential operators on tori 
and prove that up to a constant multiple, it gives the unique 
\emph{continuous} trace on the algebra of classical  
pseudodifferential operators. Let us denote the space of all classical 
Euclidean pseudodifferential 
operators with integral orders on $\mathbb{T}^n$ by 
$\Psi^{\infty}_{{\rm cl}}
(\mathbb{T}^n \times \mathbb{R}^n)$. A linear functional 
$\varphi: \Psi^{\infty}_{{\rm cl}}(\mathbb{T}^n \times \mathbb{R}^n) \to 
\mathbb{C}$ is said to be \emph{continuous} if there exists an integer $N$
such that $\varphi$ vanishes on any pseudodifferential operator of order 
less than $N$.
 
\newtheorem{first}{Definition}[section]
\begin{first} \label{first}
Let $\sigma \in S^m_{{\rm cl}}(\mathbb{T}^n \times \mathbb{R}^n)$ 
be such 
that it has  an 
asymptotic expansion of the form 
\begin{equation}
\sigma(x,\xi)\sim \sum_{j=0}^{\infty} \sigma_{m-j}(x,\xi) \,\,\,  
\textnormal{as} \,\,\, \xi \to \infty, 
\end{equation}
with positively homogeneous terms as 
in \eqref{classicaldef}. Then the noncommutative
residue of the pseudodifferential operator $T_{\sigma}$ is defined  by
\begin{equation} \label{ncresiduedef}
\textnormal{Res}\,(T_{\sigma}) = 
\int_{\mathbb{S}^{n-1}}\int_{\mathbb{T}^n} 
\sigma_{-n}(x, \xi) \, \textnormal{d}x \, \textnormal{d}\Omega,
\end{equation}
where $\mathbb{S}^{n-1}$ is the unit sphere in $\mathbb{R}^n$ centered at 
the origin, and 
$\textnormal{d} \Omega$ is the usual surface measure on the sphere.
\end{first}

Clearly, the noncommutative residue 
$\textnormal{Res}: \Psi^{\infty}_{{\rm cl}}(\mathbb{T}^n \times 
\mathbb{R}^n) 
\to \mathbb{C}$ is a linear functional that vanishes 
on any operator of order less than $-n$. In particular, it vanishes on the 
\emph{infinitely smoothing operators} 
$\Psi^{-\infty}_{{\rm cl}}(\mathbb{T}^n \times \mathbb{R}^n)$ defined by   
$$\Psi_{{\rm cl}}^{-\infty}(\mathbb{T}^n\times \mathbb{R}^n)=\cap_{m \in 
\mathbb{Z}} \Psi^{m}_{cl}(\mathbb{T}^n \times \mathbb{R}^n).$$

\newtheorem{trace}[first]{Theorem}
\begin{trace}
The noncommutative residue $\textnormal{Res} : \Psi^{\infty}_{cl}(\mathbb{T}^n 
\times \mathbb{R}^n) \to \mathbb{C}$ is a trace
 on the algebra of classical Euclidean pseudodifferential
operators on $\mathbb{T}^n$. Moreover, it vanishes on the infinitely smoothing 
operators and up to multiplication by a constant, it is the 
unique continuous trace on this algebra.
\begin{proof}
Let $\sigma \in S^m_{cl}(\mathbb{T}^n \times \mathbb{R}^n)$ and let  $\tau 
\in 
S^{m'}_{cl}(\mathbb{T}^n
\times \mathbb{R}^n)$, where $m, m' \in \mathbb{Z}$. Suppose that $\sigma$ 
and $\tau$  have asymptotic expansions given by
$$
\sigma(x,\xi)\sim \sum_{j=0}^{\infty} \sigma_{m-j} (x, \xi)$$ and  
$$\tau(x, \xi) \sim \sum_{k=0}^{\infty} \tau_{m'-k}(x, \xi)$$  
as  $|\xi| \to \infty$, 
where $\sigma_{m-j}$ and $\tau_{m'-k}$ are homogeneous of order $m-j$ and 
$m'-k$ 
respectively. Using the product formula \eqref{compositionsymbol}, the 
symbol $\lambda$ of $T_{\sigma} T_{\tau}$ 
has the asymptotic expansion given by
 \begin{eqnarray}
\lambda &\sim& \sum_{\gamma} \frac{1}{\gamma !} 
(\partial_\xi^{\gamma} 
\sigma) (D_x^{\gamma} \tau) \nonumber \\
&\sim& \sum_{\gamma} \sum_{j=0}^{\infty} \sum_{k=0}^{\infty} 
\frac{1}{\gamma !} (\partial_{\xi}^\gamma 
\sigma_{m-j})( D_x^\gamma \tau_{m'-k}). \nonumber
\end{eqnarray}
Since $ \partial_{\xi}^\gamma \sigma_{m-j} D_x^\gamma \tau_{m'-k}$ is 
homogeneous of order $m-j-|\gamma|+m'-k$, by 
the definition of the noncommutative residue given by \eqref{ncresiduedef}, 
we have
\begin{eqnarray}
&&\text{Res}\,(T_\sigma T_\tau)\nonumber\\ &=& \sum_{\substack{\gamma, j, 
k \geq 0 \\ 
m-j-|
\gamma|+m'-k=-n}} \frac{1}{\gamma !} 
\int_{\mathbb{S}^{n-1}} \int_{\mathbb{T}^n}(\partial_{\xi}^\gamma 
\sigma_{m-j})(x, \xi) 
(D_x^\gamma \tau_{m'-k}) (x, \xi) \, \text{d}x \, \text{d} \Omega. 
\nonumber
\end{eqnarray}
Similarly,
\begin{eqnarray}
&&\text{Res}\,(T_\tau T_\sigma)\nonumber\\ &=& \sum_{\substack{\gamma, k, 
j 
\geq 0 \\ m'-k-|
\gamma|+m-j=-n}} \frac{1}{\gamma !} 
\int_{\mathbb{S}^{n-1}}\int_{\mathbb{T}^n}(\partial_{\xi}^\gamma 
\tau_{m'-k})
(x, \xi) 
(D_x^\gamma \sigma_{m-j}) (x, \xi) \, \text{d}x \, \text{d}\Omega. 
\nonumber
\end{eqnarray}
Therefore we have
\begin{eqnarray}
&&\text{Res}\,(T_\sigma T_\tau) - \text{Res}\,(T_\tau T_\sigma) \nonumber 
\\
&=& \sum \int_{\mathbb{S}^{n-1}} \int_{\mathbb{T}^{n}}  \frac{1}{\gamma !} 
\{(    \partial_{\xi}^\gamma \sigma_{m-j})(x, \xi) 
(D_x^\gamma \tau_{m'-k}) (x, \xi) -  \nonumber \\
&& \qquad \qquad \qquad \qquad \qquad (D_x^\gamma \sigma_{m-j}) (x, \xi) 
(\partial_{\xi}^\gamma \tau_{m'-k})(x, \xi) 
 \} \, \text{d}x \, \text{d} \Omega, \nonumber
\end{eqnarray}
where the summation is over all nonnegative integers $j, k$ and 
multi-indices $\gamma$ 
such that
\[
m-j-|\gamma|+m'-k=-n. \nonumber
\]
Each term
\[
(\partial_{\xi}^\gamma \sigma_{m-j})(x, \xi) (D_x^\gamma \tau_{m'-k}) (x, 
\xi) 
-  (D_x^\gamma \sigma_{m-j}) (x, \xi) 
(\partial_{\xi}^\gamma \tau_{m'-k})(x, \xi) 
\]
in the above integral can be written in the form 
\[ 
\sum_{\ell=0}^n ((\partial_{\xi_\ell} A_\ell)(x, \xi) +( D_{x_\ell} 
B_\ell)(x, \xi)
 ),
\]
for some smooth maps $A_\ell$ and $ B_\ell$.  Considering the order of 
homogeneity 
$$m-j-\gamma+m'-k=-n,$$ each 
$A_\ell$ is homogeneous of order $-n+1$, and each $B_\ell$ is homogeneous 
of order $-n$ in $\xi$. In view of Lemma 
1.2 of \cite{FedGolLeiSch}, we have 
\[
\int_{\mathbb{S}^{n-1}}(\partial_{\xi_\ell} A_\ell)(x, \xi) \, \text{d} 
\Omega= 0,
\]
and from integration by parts it follows that
\[
\int_{\mathbb{T}^{n}} (D_{x_\ell}B_\ell)(x, \xi) \, \text{d}x =0. 
\]
Hence we have proved that $\text{Res}$ is a trace to the effect that
\[ 
\text{Res}\,(T_\sigma T_\tau) - \text{Res}\,(T_\tau T_\sigma) = 0.
\]
As for the uniqueness, let $\varphi : \Psi_{{\rm 
cl}}^{\infty}(\mathbb{T}^n \times \mathbb{R}^n) /
\Psi_{{\rm cl}}^{-\infty}(\mathbb{T}^n \times \mathbb{R}^n)\to \mathbb{C}$ 
be a 
continuous trace on the algebra of classical pseudodifferential operators 
modulo the infinitely smoothing 
operators. For any classical symbol $\sigma$, the symbol of 
$T_{\xi_\ell} T_\sigma - T_\sigma T_{\xi_\ell}$ is 
equivalent to $D_{x_\ell} \sigma$. Therefore 
\begin{equation} 
\varphi(T_{D_{x_{\ell}} \sigma}) =0, \,\,\, \ell = 1,2, \dots, n. 
\end{equation}
Also, $\varphi$ vanishes on $T_\sigma T_{e^{ix_{\ell}}} - T_{e^{ix_{\ell}}} 
T_\sigma$, of which the symbol is equivalent to 
\[
(\partial_{\xi_{\ell}} \sigma) e^{ i x_{\ell}} + \frac{1}{2!} (\partial^2_{\xi_{\ell}} \sigma)  e^{ i x_{\ell}} + 
\frac{1}{3!} (\partial^3_{\xi_{\ell}} \sigma) e^{ i x_{\ell}} + \cdots .
\]
By iteration, namely, by using the same argument for $\partial_{\xi_\ell} 
\sigma$ instead of $\sigma$ and so on,  and using the continuity 
of $\varphi$, it follows that 
\begin{equation} \label{vanishxi}
\varphi(T_{\partial_{\xi_{\ell}} \sigma}) = 0, \,\,\, \ell = 1,2, \dots, 
n. 
\end{equation}
Now, assume that $$\sigma (x,\xi)\sim \sum_{j=0}^{\infty} \sigma_{m-j} (x, 
\xi)$$ 
as $|\xi| \to \infty$, where each $\sigma_{m-j}$ is 
homogeneous of order $m-j$. If $m-j \neq -n$, using Euler's identity, there 
are smooth map $h_{\ell, m-j}$ such that 
\[
\sigma_{m-j} = \sum_{\ell=0}^{n} \partial_{\xi_\ell} h_{\ell, m-j}.
\] 
Setting 
\[  h_\ell \sim \sum_{\substack{j \geq 0 \\ m-j \neq - n}} h_{\ell, m-j}, 
\,\,\, \ell=1,2,  \dots, n,\nonumber\]
we have 
\[ \sigma \sim \sigma_{-n} +  \sum_{\ell=1}^n \partial_{\xi_\ell} h_\ell. 
\nonumber \]
Now, from \eqref{vanishxi} it follows that 
\[\varphi(T_\sigma) = \varphi (T_{\sigma_{-n}}).\]
We write
\[ \sigma_{-n}(x, \xi) = \sigma_{-n}(x, \xi) - r(x) |\xi|^{-n} +  r(x) 
|\xi|^{-n}, \]
where 
\[ r(x) =  \frac{1}{ |\mathbb{S}^{n-1}|}
\int_{\mathbb{S}^{n-1}} \sigma_{-n}(x, \xi) \, \textnormal{d} \Omega,\] 
where $|\mathbb{S}^{n-1}|$ is the surface measure of $\mathbb{S}^{n-1}.$
Since $$\int_{\mathbb{S}^{n-1}} (\sigma_{-n}(x, \xi) -r(x) |\xi|^{-n} ) \, 
\text{d} \Omega =0,$$ it follows from Lemma 1.3 of \cite{FedGolLeiSch} 
that  $\sigma_{-n}(x, \xi) -r(x) |\xi|^{-n}$ can be written as a sum of 
partial derivatives. Hence, using \eqref{vanishxi}, $\varphi$ 
vanishes on the corresponding operator. Therefore
\[ \varphi(T_{\sigma}) = \varphi (  T_{r(x) |\xi|^{-n}}). \]
Now, we  consider the linear map $\chi : C^{\infty}(\mathbb{T}^{n}) \to 
\mathbb{C}$ 
defined by 
\[ \chi(f) = \varphi(T_{f |\xi|^{-n}}). \nonumber  \] 
Since $$\chi(\partial_{x_{\ell}}f) =0$$ for all $f\in C^{\infty}
(\mathbb{T}^{n})$ and $\ell=1,2, \dots, n$, there exists a constant $c$ 
such that 
\[ \chi(f) = c \int_{\mathbb{T}^n} f(x) \, \textnormal{d}x, \,\,\,  f 
\in C^{\infty}(\mathbb{T}^n). \nonumber \]
Therefore 
\begin{eqnarray}
\varphi(T_{\sigma}) &=&    \int_{\mathbb{T}^n} r(x) \, \textnormal{d}x 
\nonumber \\
&=& \frac{c}{|\mathbb{S}^{n-1}|} \int_{\mathbb{T}^n}
\int_{\mathbb{S}^{n-1}} \sigma_{-n}(x, \xi) 
\, \textnormal{d} \Omega \, \textnormal{d}x \nonumber \\
&=&    \frac{c}{|\mathbb{S}^{n-1}|} \textnormal{Res}\,
(T_{\sigma}). \nonumber 
\end{eqnarray}
Therefore $\varphi$ is a constant multiple of the noncommutative residue.
\end{proof}
\end{trace}

\section{The Noncommutative Two Torus}
\setcounter{equation}{0}
In this section we first recall Connes' pseudodifferential calculus for the 
canonical dynamical system associated to the 
noncommutative two torus $A_\theta$, $\theta \in \mathbb{R}$ \cite{con}. Then we  show 
that in the  case $\theta =0$, 
the noncommutative residue for classical pseudodifferential operators on 
$A_\theta$ defined in \cite{FatWon} coincides 
with the noncommutative residue defined in Section \ref{residueeuclidean}.

By definition, for a fixed $\theta \in \mathbb{R}$, $A_{\theta}$ is the 
universal unital $C^*$-algebra  generated by two 
unitaries $U$ and $V$ satisfying
\[VU=e^{2 \pi i \theta} UV.\]
There is a continuous action of $\mathbb{T}^2$, $\mathbb{T}= \mathbb{R}/2\pi 
\mathbb{Z}$, on $A_{\theta}$ by $C^*$-algebra
automorphisms  $\{ \alpha_s\}$, $s\in \mathbb{R}^2$, defined by
\begin{equation} \label{dynamic}
\alpha_s(U^mV^n)=e^{is\cdot(m,n)}U^mV^n.
\end{equation}
The space of smooth elements for this action, $i.e.,$  elements $a 
\in A_{\theta}$ for which the map $s \mapsto
\alpha_s (a)$  is $C^{\infty}$, is  denoted by
$A_{\theta}^{\infty}$. It is a dense subalgebra of $A_{\theta}$  which can be 
alternatively
described as the algebra of elements in $A_{\theta}$
whose (noncommutative) Fourier expansion has rapidly decreasing 
coefficients. More precisely,
\[A_{\theta}^{\infty}= \left\{\sum_{m,n\in \mathbb{Z}}a_{m,n}U^mV^n:
\sup_{m,n \in \mathbb{Z}} (|m|^k|n|^q|a_{m,n}|)< \infty, \,  k, q \in 
\mathbb{Z}  \right \}.\]

There exist two derivations $\delta_1:A_\theta^\infty\to A_\theta^\infty$ 
and $\delta_2:A_\theta^\infty\to A_\theta^\infty$ corresponding to the 
above action of $\mathbb{T}^2$ on $A_\theta$, which are, respectively, the 
analogs of the differential operators $-i\frac{\partial}{\partial x}$ and 
$-i\frac{\partial}{\partial y}$ on smooth functions on $\mathbb{T}^2.$ 
These derivations are fixed by
$$\delta_1(U)=U,\quad \delta_1(V)=0,$$
and
$$\delta_2(U)=0,\quad \delta_2(V)=V.$$
Moreover, for $j=1,2,$ we have
$$\delta_j(a^*)=-\delta_j(a)^*$$
for all $a\in A_\theta^\infty.$

There is a normalized trace $\mathfrak{t}$ on $A_{\theta}$ that turns 
out to be positive and faithful. This means that $$\mathfrak{t}\,(1)=1$$ and 
$${\mathfrak{t}}\,(a^*a) > 
0$$ 
for all nonzero $a \in A_{\theta}$. The restriction of ${\mathfrak{t}}$ to 
$A_\theta^\infty$ is given by
$${\mathfrak{t}}\left(\sum_{j,k\in 
\mathbb{Z}}a_{j,k}U^jV^k\right)=a_{0,0}$$ for all 
$$\sum_{j,k\in \mathbb{Z}}a_{j,k}U^jV^k\in A_\theta^\infty.$$

For any integer $n$, a smooth map $\sigma: \mathbb{R}^2 \to 
A_{\theta}^{\infty}$ 
is  said
to be a symbol of order $n$ \cite{con}, if for all nonnegative integers $i_1, i_2, j_1,
j_2,$ there exists a positive constant $C$, depending on $i_1$ $i_2$, 
$j_1$ and $j_2$ only, such that
\[ ||\delta_1^{i_1} \delta_2^{i_2} ((\partial_1^{j_1} \partial_2^{j_2} 
\sigma)(\xi) 
||
\leq c (1+|\xi|)^{n-j_1-j_2},\]
 and if there exists a smooth map $k: \mathbb{R}^2 \to
A_{\theta}^{\infty}$ such that
\[\lim_{\lambda \to \infty} \lambda^{-n} \sigma(\lambda\xi_1, 
\lambda\xi_2) = 
k (\xi_1, \xi_2)\] for all $(\xi_1,\xi_2)$ in $\mathbb{R}^2.$
The space of symbols of order $n$ is denoted by $S^n$.

To a symbol $\sigma$ of order $n$, we associate an operator on 
$A_{\theta}^{\infty}$ \cite{con},
denoted by $T_{\sigma,\theta}$ and  given by
\begin{equation} \label{NCtorus} 
T_{\sigma,\theta}(a) = (2 \pi)^{-2} \int_{\mathbb{R}^2} 
\int_{\mathbb{R}^2} e^{-is 
\cdot \xi} \sigma(\xi) \alpha_s(a) 
\,\text{d}s \,
\text{d}\xi,   \,\,\,\,\,\,\, a \in A_{\theta}^{\infty},
\end{equation}
where $\alpha_s$ is given by \eqref{dynamic}.

A symbol $\sigma$ in $S^n$ is 
said to be a {\it classical} symbol and we write $\sigma \in S^n_{\rm cl}$ if
it admits an asymptotic expansion of the form
$$\sigma(\xi) \sim \sum_{j=0}^{\infty} \sigma_{n-j}(\xi)$$  as $|\xi|\to 
\infty$,
where for  each $j=0,1,2, \dots$, $\sigma_{n-j} : \mathbb{R}^2 \setminus 
\{ 
0 
\} \to A_{\theta}^{\infty}$ is
smooth and positively homogeneous of order $n-j$.
Then we define the
noncommutative residue ${\rm Res}\,(T_{\sigma,\theta})$  of 
$T_{\sigma,\theta}$ by
\begin{equation} \label{resdef}
{\rm Res}\,(T_{\sigma,\theta}) = \int_{\mathbb{S}^1} 
{\mathfrak{t}}\,(\sigma_{-2}(\xi))\,\textnormal{d}\Omega,
\end{equation}
where $\textnormal{d}\Omega$ is the Lebesgue measure on the unit circle 
${\mathbb{S}}^1$ centered at the origin \cite{FatWon}. The space of pseudodifferential operators on 
$A_\theta$ form an algebra \cite{con}, and it is shown in \cite{FatWon} that the above noncommutative 
residue is the unique continuous trace on the classical operators.

Now, let us consider the noncommutative two torus $A_{\theta}$ for the 
case 
$\theta=0$. By definition, we can assume that $A_0$
is the $C^*$-algebra generated by the functions $U$ and $V$ on 
$\mathbb{R}^2$ defined by
$$
U(x,y) = e^{ix}$$ and $$ V(x,y)= e^{iy}$$ for all $(x,y)\in \mathbb{R}^2.$ 
Since these functions are $2 \pi$-periodic in both variables, we can consider 
them as smooth functions defined on the 
two torus $\mathbb{T}^2$, $\mathbb{T}= \mathbb{R}/2\pi \mathbb{Z}$. For 
all 
smooth function $f$ defined on $\mathbb{T}^2$, 
we use its Fourier expansion to obtain
$$
f = \sum_{m, n \in \mathbb{Z}} a_{m, n} U^m V^n, $$ 
where 
$$
a_{m, n} = \frac{1}{4 \pi^2} \int_{0}^{2 \pi} \int_{0}^{2 \pi} 
f(x, y) e^{-imx} e^{-iny} \, \textnormal{d}x\, \textnormal{d}y.$$
Hence in the case $\theta=0$, $A_\theta$ is the algebra of continuous 
functions on the ordinary two torus and $A_{\theta}^\infty$ 
is the algebra of smooth functions on $\mathbb{T}^2$.

So, if we translate the function $f$ by $s=(s_1, s_2)\in \mathbb{R}^2$, 
and denote the result by  $T_sf$, then
$$
T_s f = \sum_{m,n\in \mathbb{Z}} e^{i s\cdot (m, n)} a_{m, n} U^m V^n $$
because we have
\begin{eqnarray} 
&&\int_{0}^{2 \pi} \int_{0}^{2 \pi} f(x+s_1, y+s_2) e^{-imx} e^{-iny} 
\textnormal{d}x\,\textnormal{d}y \nonumber \\
&=&  e^{ims_1} e^{ins_2} \int_{0}^{2 \pi} \int_{0}^{2 \pi} f(x, y) e^{-imx} 
e^{-iny} \textnormal{d}x\,\textnormal{d}y. \nonumber
\end{eqnarray}
Thus, in the case $\theta=0$, the action $\alpha_s$ on $A_{\theta}$ 
described in \eqref{dynamic} is 
just the translation of functions by $s$. Now by a simple change of variable, namely by passing to $s=y-x$, 
we can observe the following identity for a variant of formula \eqref{pseudodef} for $n=2$, and formula \eqref{NCtorus}. 
In fact, we have
\begin{eqnarray}
\int_{\mathbb{R}^2} \int_{\mathbb{R}^2} e^{ i (x-y) \cdot \xi} \sigma(x, \xi) f(y) \, \text{d}y \, \text{d}\xi & = & 
\int_{\mathbb{R}^2} \int_{\mathbb{R}^2} e^{ -i s\cdot \xi} \sigma(x, \xi) f(x+s) \, \text{d}s \, \text{d}\xi  \nonumber \\
&=&\int_{\mathbb{R}^2} \int_{\mathbb{R}^2} e^{ -i s\cdot \xi} \sigma(x, \xi) \alpha_s(f)(x) \, \text{d}s \, \text{d}\xi. \nonumber
\end{eqnarray}
Moreover, in the commutative case $\theta=0$, the 
trace $\mathfrak{t}$ on $A_0$ amounts to 
integration of continuous functions on the ordinary two torus 
$\mathbb{T}^2.$ 
We also note that in this case,
the $C^*$-algebra norm is given by the supremum norm of 
continuous 
functions on $\mathbb{T}^2$. Therefore the definition of the symbols 
given above coincide with the one given in Section \ref{residueeuclidean}. 
Note that here we are identifying the complex-valued smooth functions defined on 
$\mathbb{T}^2 \times \mathbb{R}^2$ with the smooth functions from $\mathbb{R}^2$ to $A_\theta^\infty$.

Considering the above observations and the fact  that the noncommutative residues \eqref{ncresiduedef} 
and \eqref{resdef} are defined on algebras of classical pseudodifferential symbols with multiplications 
induced from 
composition of pseudodifferential operators, it is clear that for $n=2$, the noncommutative 
residue \eqref{ncresiduedef} is the semiclassical limit of \eqref{resdef} defined on the pseudodifferential 
symbols on the noncommutative two torus 
$A_\theta$.  We record this result in the following.

\newtheorem{coincide}[first]{Theorem}
\begin{coincide}
In the case $\theta=0$, the noncommutative residue defined by \eqref{resdef} 
on classical pseudodifferential operators on 
the noncommutative two torus $A_{\theta}$, coincides with the noncommutative 
residue defined by \eqref{ncresiduedef} on the classical Euclidean pseudodifferential 
operators on $\mathbb{T}^2$.  
\end{coincide}

\noindent
Department of Mathematics and Statistics \\York University \\ Toronto, Ontario, Canada, M3J 1P3 \\
ffathiza@mathstat.yorku.ca

\begin{thebibliography}{9}

\bibitem{Adl} M. Adler, On a trace functional for formal 
pseudo-differential operators and the symplectic structure of Korteweg--de 
Vries type equations, {\it Invent. Math.} {\bf 50} (1978/79), 219--248.


\bibitem{con} A. Connes,  {$C^*$-alg\`ebres et g\'eom\'etrie 
diff\'erentielle}, {\it C.R. Acad. Sc. Paris
S\'er. A-B} {\bf 290} (1980), A599--A604.

\bibitem{conmos} A. Connes and H.  Moscovici,  { The local index formula in
noncommutative geometry.}
{\it  Geom. Funct. Anal. 5}, no. 2 (1995), 174--243.


\bibitem{FatWon} F. Fathizadeh and M. W. Wong, Noncommutative residues for pseudo-differential 
operators on the noncommutative two-torus,  \emph{Journal of
 Pseudo-Differential Operators and Applications}, {\bf 2}(3)
 (2011), 289--302.

\bibitem{FedGolLeiSch} B. Fedosov, F. Golse, E. Leichtnam and E. Schrohe, 
{The noncommutative residue for manifolds
with boundary}, {\it J. Funct. Anal.} {\bf 142}(1) (1996), 
1--–31.



\bibitem{Man} Ju. I. Manin, Algebraic aspects of nonlinear differential 
equations, {\it J. Soviet Math.} {\bf 11} (1979), 1--122. 


\bibitem{RuzTur} M. Ruzhansky and V. Turunen, Pseudo-Differential Operators 
and Symmetries, Birkh\"auser, 2010.

\bibitem{Wod}
M. Wodzicki, {Noncommutative residue. I. Fundamentals,}
in {\it $K$-theory, Arithmetic and Geometry (Moscow, 1984--1986)},  
Lecture Notes in Mathematics {\bf 1289}, Springer, 1987, 320--399.


\end{thebibliography}
\end{document}